\documentclass[11pt]{amsart}
\usepackage[utf8]{inputenc}
\usepackage{amssymb}
\usepackage{amscd}
\usepackage{curves}
\usepackage{epsfig}
\usepackage[pass]{geometry}
\usepackage{graphicx}
\usepackage{pst-plot}
\usepackage{pst-pdf}
\usepackage{verbatim}
\usepackage{mathrsfs}
\usepackage{amsfonts}
\usepackage{tikz}
\usepackage{bbm}
\usetikzlibrary{matrix}
\usepackage{latexsym,eucal,amsfonts,amssymb,amsmath,graphicx}
\numberwithin{equation}{section}
\newtheorem{theorem}{Theorem}[section]

\newtheorem{cor}[theorem]{Corollary}

\def \bpf {\begin{proof}}
\def \epf {\end{proof}}
\def \beq {\begin{equation*}}
\def \eeq {\end{equation*}}
\def \bsp{\begin{split}}
\def \esp{\end{split}}

\def \ha {\frac{1}{2}}

\def \mca {{\mathscr A}}
\def \mcb {{\mathscr B}}
\def \mcc {{\mathscr{C}}}

\def \mce {{\mathscr E}}

\def \mcl {{\mathscr L}}

\def \mcn {{\mathscr N}}

\def \mcp {{\mathscr P}}

\def \mcu {{\mathscr U}}

\def \mbr {{\mathbb R}}
\def \mbs {{\mathbb S}}

\def \supp {\text{supp }}

\def \eps {\epsilon}   
   
\def \La {\Lambda}

\def \lap {\Delta}
\def \p {\partial}

\def \eps {\epsilon}

\def \ha {\frac{1}{2}}

\def \beqq {\begin{equation}}
\def \eeqq {\end{equation}}

\def \sym {\text{Sym}}
\def \WF {\text{WF}}

\begin{document}
\title{Parametrices for the Light Ray Transform on Minkowski Spacetime}
\author{Yiran Wang}
\address{Yiran Wang
\newline
\indent Department of Mathematics, University of Washington,
\newline
\indent Box 354350, Seattle, WA 98195-4350 
\newline
\indent and 
\newline 
 \indent Institute for Advanced Study, the Hong Kong University of Science and Technology
 \newline
 \indent Lo Ka Chung Building, Lee Shau Kee Campus, Clear Water Bay, Kowloon, Hong Kong}
\email{wangy257@math.washington.edu}
\begin{abstract} 
We consider restricted light ray transforms arising from an inverse problem of finding cosmic strings. We construct a relative left parametrix for the transform on two tensors, which recovers the space-like and some light-like singularities of the two tensor. 
\end{abstract}
\maketitle

\section{Introduction}
Let $(M, g)$ be a smooth Lorentzian manifold.  A smooth curve $\gamma: \mbr \rightarrow M$ is called a light ray if $g(\dot \gamma(t), \dot \gamma(t)) = 0, t\in \mbr$, where $\dot \gamma$ denotes the covariant derivative along $\gamma$. Let $\mcc$ be the set of light rays on $M$ and $\sym^2$ denote the bundle of symmetric $2$-tensors on $M$. For $f\in C_0^\infty(M; \sym^2)$, we consider the light ray transform 
\beqq\label{eqlf}
L_{\mcc} f(\gamma) = \int_\mbr f_{lm}(\gamma(t)) \dot \gamma^l(t) \dot \gamma^m(t) dt, \ \ \gamma \in \mcc.
\eeqq
Hereafter, the Einstein summation convention is used i.e.\ summation is over repeated indices. On $2+1$ dimensional Minkowski space-time, this transform was studied by Guillemin \cite{Gu} and it was envisioned to have applications in \lq\lq{}cosmological X-ray tomography\rq\rq{}, see the  concluding remarks  \cite[Section 17]{Gu}. Recently in \cite{LOSU}, such transforms naturally arise from an inverse problem of detecting singularities of  the Lorentzian metric of the Universe using Cosmic Microwave Background (CMB) radiation measurements. In particular, let $(M, g)$ be a Friedmann-Lema\^itre-Robertson-Walker (FLRW) type model for the Universe. For a small parameter $\eps$, consider a family of metrics on $M$: \beq
g_\eps = g + \eps f + O(\eps^2)\eeq
 representing small perturbations of $g$. In \cite{LOSU}, it is demonstrated that one can obtain a restricted light ray transform of $f$ from the linearization of the CMB measurements. Then it is proved in Theorem 4.4 that one can recover the space-like singularities of $f$. However, as already noted in \cite{LOSU}, light-like singularities are of great interest as they correspond  to gravitational waves which may be caused for example by cosmic strings. We address this problem in this note.

For restricted geodesic ray transforms on functions (the light ray transform being an example), there is a microlocal framework developed by Greenleaf-Uhlmann \cite{GrU1, GrU2', GrU90, GrU2, GrU3} to understand their mapping properties. We combine it with some calculations in \cite{LOSU} to show that the normal operator of the light ray transform is a paired Lagrangian distribution and construct parametrices on the elliptic part. This allows us to obtain a relative left parametrix for the restricted light ray transform, which among other things recovers the space-like and some light-like singularities of the metric perturbations. We remark that time-like singularities are in the kernel of $L_{\mcc}$ and there is a good physical explanation for one not being able to determine them from the light ray transform, see \cite[Page 188]{Gu}, \cite{LOSU} and the interesting work of Stefanov \cite{St} on support theorems of the light ray transform.

The paper is organized as follows. In Section \ref{main}, we state the main results after setting up the problem. We show in Section \ref{secnorm} that the normal operator is a paired Lagrangian distribution and we construct the parametrix in Section \ref{secpara}.\\

\textbf{Acknowledgement:} The author sincerely thanks Prof.\ Gunther Uhlmann for suggesting the problem and for many helpful discussions. He is also grateful to Prof.\ Allan Greenleaf for reference \cite{GrS} and related comments.

\section{The main results}\label{main}
\begin{figure}[htbp]
 \centering
\includegraphics[scale=0.55]{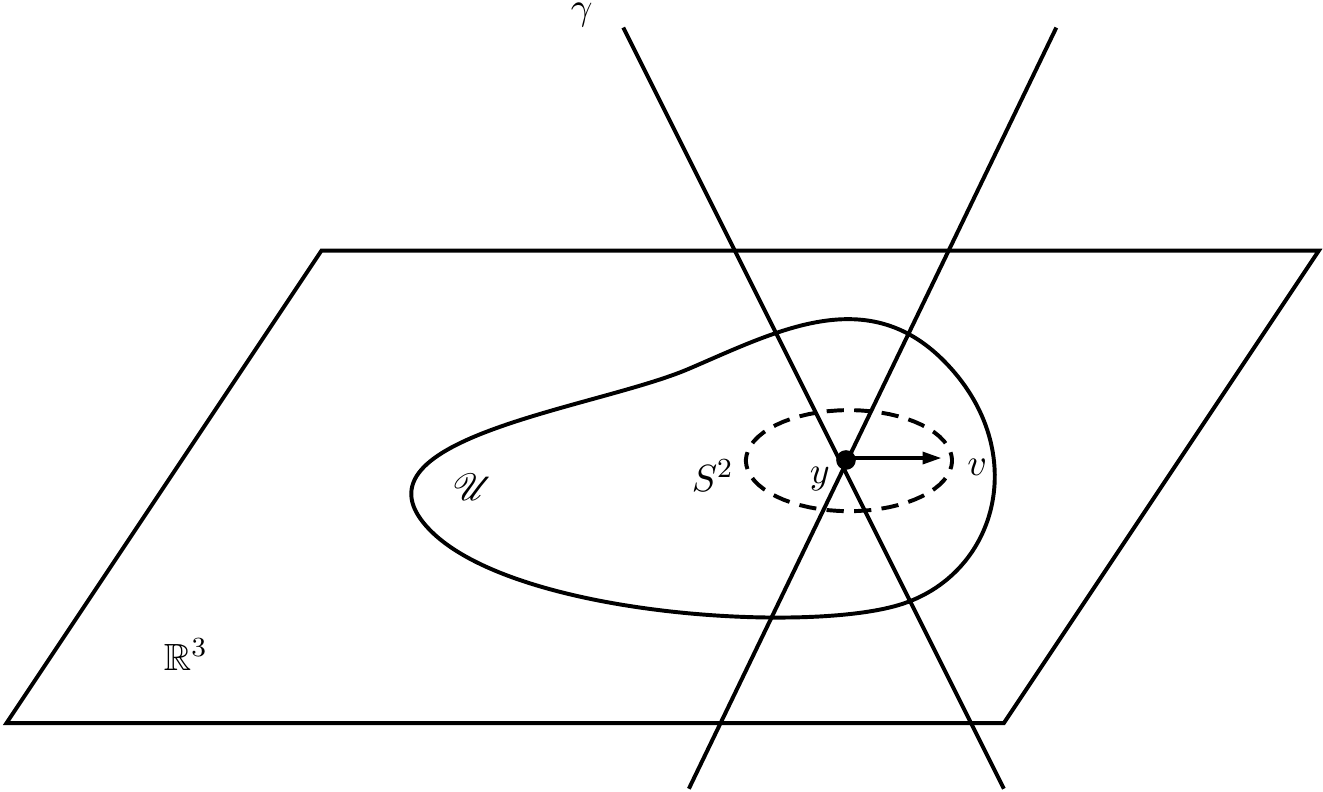}
\caption{Illustration of complex $\mcc_0$}
\label{fig1}
\end{figure}

It is known that a FLRW type space-time is conformal to the Minkowski space-time. Since conformal diffeomorphisms preserve light-like geodesics, as discussed in \cite{LOSU}, it suffices to consider light ray transforms on  
\beq
M = \mbr^{3+1}, \ \ g = -dt^2 + \sum_{i = 1}^3(dx^i)^2,
\eeq
where $x = (t, x\rq{}) = (x^0, x^1, x^2, x^3)$ denotes the coordinates on $M$. In this case, the light rays are straight lines and we denote by $\mcc$ the set of light rays. As demonstrated in \cite[Lemma 4.3]{LOSU}, the light ray transform $L_{\mcc}$ defined as \eqref{eqlf} has a non-trivial null space given by 
\beq
\mcn = \{cg + d^s \omega; c\in \mce\rq{}(M), w \in \mce\rq{}(M; \La^1)\},
\eeq
where $\mce\rq{}(M)$ denotes the space of distributions with compact support, $d^s$ is the symmetric differential given in local coordinates by 
\beq
(d^s\omega)_{ij} = \frac{(\nabla_i \omega)_j + (\nabla_j \omega)_i}{2}, \ \ i, j = 0, 1, 2, 3,
\eeq
with $\nabla_i$  the covariant derivative, and $\La^1$ denotes the bundle of one forms. 
Let $\mcu$ be an open set of $\mbr^3$ and define the line complex 
\beq
\mcc_0 = \{\gamma \in \mcc: \gamma \cap \mcu \neq \emptyset  \}
\eeq
i.e.\ collection of all light rays intersecting $\mcu$, see Figure \ref{fig1}. We denote by $L_{\mcc_0} = L_{\mcc}|_{\mcc_0}$ the restricted light ray transform on $\mcc_0$. To describe the null space of $L_{\mcc_0}$, we denote $\mcl(\mcu) = \{p \in M: \text{there exists $q\in \mcu$ and a light ray joining $p$ and $q$}\}$. Then we observe that $L_{\mcc_0}f = 0$ if $f$ is supported on $M\backslash \mcl(\mcu)$. Also, $L_{\mcc_0}f = L_{\mcc} f$ for $f$ supported in $\mcl(\mcu)$. For any $\mcb \subset \mbr^3$ we denote
\beq
\mcp(\mcb) = \{f \in \mce'(M; \sym^2)\backslash \mcn: \supp (f) \subset \mcl(\mcb)\},
\eeq
then we observe that $L_{\mcc_0}$ is injective on $\mcp(\mcu)$.

The microlocal nature of $L_{\mcc_0}$ is well-understood. Let 
\beq
Z = \{(\gamma, x)\in \mcc_0\times M: x\in \gamma\}
\eeq
be the point-line relation. Then the Schwartz kernel of $L_{\mcc_0}$ is $\delta_Z$ the delta distribution on $\mcc_0 \times M$ supported on $Z$. Hence we know from H\"omander's theory  that $L_{\mcc_0}$ is a Fourier integral operator of order $-3/4$ associated with the canonical relation $N^*Z'$ (see \eqref{eqcan}). Although we do not explore this point here, the operator should fit into the framework in \cite{GrU3}, see also \cite{GrS}. In Section \ref{secnorm}, we use a more direct approach to show that the Schwartz kernel of the normal operator $L_{\mcc_0}^t\circ L_{\mcc_0}$ is a paired Lagrangian distribution and we obtain the Sobolev estimate of $L_{\mcc_0}$, see Theorem \ref{main0}.

To state the main result, we need to describe the two Lagrangians associated to the normal operator. Let $T^*M$ be the cotangent bundle and $(x, \xi)$ be the coordinate for $T^*M$ where $\xi = (\xi_0, \xi_1, \xi_2, \xi_3)$. Consider $p(x, \xi) = g(\xi, \xi) = -\xi_0^2 + \xi_1^2 + \xi_2^2 + \xi_3^2$ the (dual) metric function on $T^*M$.  We denote $L^*M = \{(x, \xi)\in T^*M: g(\xi, \xi) = 0\}$ the light-like covectors, $\Omega^+ = \{(x, \xi)\in T^*M: g(\xi, \xi)> 0 \}$ the space-like covectors and $\Omega^- = \{(x, \xi)\in T^*M:   g(\xi, \xi)<0 \}$ the time-like covectors.  
Then we can decompose $T^*M = \Omega^+ \cup \Omega^- \cup L^*M.$    Let $\omega = d\xi\wedge dx$ be the canonical two form on $T^*M$. The Hamilton vector field of $p$ denoted by $H_p$ is defined through 
\beq
H_p = \sum_{i = 0}^3 \frac{\p p}{\p \xi_i} \frac{\p}{\p x^i} -  \frac{\p p}{\p x^i} \frac{\p}{\p \xi^i} = 2(-\xi_0 \frac{\p}{\p x^0} + \sum_{i = 1}^3\xi_i \frac{\p}{\p x^i}).
\eeq
The integral curves of $H_p$ in  $L^*M$ are called null bicharacteristics. It is well known that their projections to $M$ are light-like geodesics. We denote $\lap = \{(x, \xi, x, -\xi) \in T^*M\backslash 0\times T^*M\backslash 0\}$ and $\Sigma = \{(x, \xi, x, -\xi) \in L^*M\backslash 0\times L^*M\backslash 0\}$ where $0$ stands for the zero section. We let $\La$ be the flow out of $\Sigma$ meaning 
\beq
\La = \{(x, \xi, y, -\eta) \in (L^*M\backslash 0) \times (L^*M\backslash 0): (x, \xi) = \exp t H_p(y, \eta), \text{ for some $t\in \mbr$}\}. 
\eeq
Then $\lap$  and $\La$ are Lagrangian subamanifolds of $T^*(M\times M)$ and they form a pair of cleanly intersecting Lagrangians in the following sense: two Lagrangians $\La_0, \La_1 \subset T^*M\backslash 0$ intersect cleanly  if 
\beq
T_p\La_0\cap T_p\La_1  = T_p(\La_0\cap \La_1),\ \ \forall p\in \La_0\cap \La_1.
\eeq

Now we briefly recall the notion of Lagrangian and paired Lagrangian distributions. Let $\La$ be a smooth conic Lagrangian submanifold of $T^*M\backslash 0$. We denote by $I^\mu(M; \La)$ the space of Lagrangian distributions  of order $\mu$ on $M$ associated with $\La$.  For two Lagrangians $\La_0, \La_1 \subset T^*M\backslash 0$ intersecting cleanly at a codimension $k$ submanifold, the space of paired Lagrangian distributions associated with $(\La_0, \La_1)$ is denoted by $I^{p, l}(M; \La_0, \La_1)$, We use $I^{p, l}(\La_0, \La_1)$ when the background manifold is clear. By abuse of notations, we also use $I^{p, l}(\La_0, \La_1)$ for section valued distributions in $\sym^2$. We know (from e.g.\ Prop.\ 3.1 of \cite{GrU1}) that if $u \in I^{p, l}(\La_0, \La_1)$, then $u\in I^{p+l}(\La_0\backslash \La_1)$ and $u \in I^p(\La_1\backslash \La_0)$. So $u$ has well-defined symbols on each Lagrangian. 
 
For any subset $\mca$ of $T^*M$, we let $\mathbbm{1}_{\mca}$ be the microlocal cut-off  defined as
\beq
\mathbbm{1}_{\mca}f(x) = \frac{1}{(2\pi)^4} \int_{\mbr^4}\int_{\mbr^4}e^{i(x-z)\eta}\chi_{\mca}(x, \eta) f(z)dzd\eta, 
\eeq
where $\chi_\mca$ is the characteristic function for $\mca$ and $f\in \mce'(M; \sym^2)$. Our main result is 

\begin{theorem}\label{main1}
There exists a relative left parametrix $A$ for $L_{\mcc_0}$ such that
\beq
A\circ L_{\mcc_0} = \mathbbm{1}_{\Omega^+}+B \text{ mod } C^\infty \text{ on } \mcp(\mcu),
\eeq
where $A = \tilde A \circ L_{\mcc}^t$ , $\tilde A \in I^{\ha, \ha}(\lap, \La)$ and $B\in I^{-\ha}(\La)$.
\end{theorem}

Using this result as a reconstruction formula and wave front analysis, we see that for $f\in  \mcp(\mcu)$,  we can recover the singularities in $f$ on space-like directions and on some light-like directions. One may not be able to recover all light-like singularities due to the error term, see a related example in \cite[Section 2]{GrU2}. However, $Bf$ contains singularities on the flow out which can be regarded as artifacts in the reconstruction. As we already mentioned, light-like singularities corresponds to gravitational waves and the artifacts may help us to identify these singularities. Furthermore, we notice that $B$ is an Fourier integral operator associated with the canonical relation $\La'$. The rank of the projection of $\La'$ to $T^*M$ drops by $1$. From H\"ormander's result on $L^2$ boundedness of Fourier integral operators  \cite[Theorem 4.3.2]{Ho0}, we conclude that if $f\in H^s(M; \sym^2)\cap \mcp(\mcu)$, then $Bf\in H^{s}(M)$. So the artifacts have the same order of Sobolev regularity as $f$ does. In a different context \cite{PUW}, the problem of reducing and enhancing the artifacts due to a similar mechanism is studied. The same strategy should work here as well. 

Away from the light-like directions, we state Theorem \ref{main1} as a corollary in the same spirit as \cite[Theorem 4.4]{LOSU}.  
\begin{cor}
For $f\in \mcp(\mcu)$ with $\WF(f)\subset \Omega^+$ and $A$ defined in Theorem \ref{main1}, we have
\beq
A\circ L_{\mcc_0}f = f \text{ mod } C^\infty \text{ in } \mcp(\mcu).
\eeq
\end{cor}

To conclude this section, we briefly review the local representations of paired Lagrangians needed for our analysis. Let's consider the space $I^{p, l}(M\times M; \lap, \La)$. By Prop.\ 2.1 of \cite{GU}, it is convenient to consider the distributions on the following model pair, which can be found in \cite{GrU90, GU, DUV}. On $\mbr^{2n} = \mbr^n_x\times \mbr^n_y$, $n = 4$, we let: $\tilde \La_0 = \{(x, \xi, x, -\xi) \in T^*\mbr^n\backslash 0 \times T^*\mbr^n\backslash 0\}$ and $\tilde \La_1 = \{(x, \xi, y, \eta): x'' = y'', \xi' = \eta' = 0, \xi'' = \eta''\}$ where $x = (x', x'')\in \mbr^k\times \mbr^{n-k}, k = 1$. So $\tilde \La_0$ intersects $\tilde \La_1$ at a codimension $k$ submanifold. In this case, we can write $u\in I^{p, l}(\tilde \La_0, \tilde \La_1)$ as 
\beq
u(x, y) = \int e^{i(x'-y'-s)\cdot \xi' + (x'' - y'')\cdot \xi'' + s\cdot \sigma} a(x, y, s, \xi, \sigma) d\sigma ds d\xi,
\eeq
with $a\in S^{p-n/2+k/2, l-k/2}(\mbr^{2n+ k}; \mbr^n, \mbr^k)$, which by definition means that for any compact set $K\subset \mbr^{2n+k}$ and multi-indices $\alpha, \beta, \gamma \geq 0$, we have
\beq
|\p_\xi^\alpha \p_\sigma^\beta \p_s^\gamma \p_x^\delta \p_y^\eps a(x, y, s, \xi, \sigma)| \leq C_{\alpha\beta\gamma\delta\eps K} (1 + |\xi|)^{p-n/2+k/2}(1 + |\sigma|)^{l-k/2}, \ \ C_{\alpha\beta\gamma\delta\eps K}  > 0.
\eeq
On $\tilde \La_0\backslash \tilde\La_1$, we can write 
\beq
u(x, y) = \int e^{i(x-y)\cdot \xi} a(x, y, 0, \xi, \xi') d\xi
\eeq
modulo a pseudo-differential operator of lower orders and the symbol $a$ is singular at $\tilde \La_0\cap \tilde\La_1$.

\section{The normal operator}\label{secnorm}
We choose a parametrization of $\mcc_0$ and find  the normal operator in the parametrization. Some of these are done in \cite{LOSU}. Let $y\in \mcu \subset \mbr^3, v\in \mbs^2 \doteq \{z\in \mbr^3: |z| = 1\}$. We let $\theta = (1, v)$ so that $\theta$ is a (future pointing) light like vector, see Fig.\ \ref{fig1}. For $\gamma \in \mcc_0$ and $\gamma\cap \mcu = (0, y)$, we can write 
\beq
\gamma(s) =   (s, y + sv), \ \ s\in \mbr.
\eeq
Then for $f\in C_0^\infty(M; \sym^2)$, we have
\beq
L_{\mcc_0}f(y, v) = \int_\mbr f_{lm}(s, y+ sv) \theta^l \theta^m ds, \ \ y\in \mcu, v\in \mbs^2.
\eeq
The point-line relation is parametrized by
\beq
Z = \{(y, v, x) \in \mcu \times S^2 \times M:  x' = y+ x^0 v\}.
\eeq
Therefore, we can find the conormal bundle $N^*Z$ and the canonical relation $C = N^*Z'$ as
\beqq\label{eqcan}
\begin{gathered}
C = \{(y, v, \eta, w, x^0, x', \xi_0, \xi') \in (T^*\mcc_0 \times T^*M)\backslash 0: y = x'-x^0v,\ \ \eta = \xi', \\
w = x^0\xi'|_{T_vS^2},\ \ \xi_0 = -\xi' v, \ \ y\in \mcu, v\in S^2, \eta \in \mbr^3, x\in M\},
\end{gathered}
\eeqq
see $(39)$ of \cite{LOSU}. Now let's consider the double fibration picture
\begin{center}
\begin{tikzpicture}
  \matrix (m) [matrix of math nodes,row sep=1.5em,column sep=1em,minimum width=1em]
  {
      & C & \\
     T^*M &  & T^*\mcc_0\\};
  \path[-stealth]
    (m-1-2) edge node [left] {$\pi$} (m-2-1)
    (m-1-2) edge node [right] {$\rho$} (m-2-3);
\end{tikzpicture}
\end{center}
If $\rho$ is an injective immersion, the double fibration satisfies the Bolker condition. In this case, the composition $L_{\mcc_0}^t\circ L_{\mcc_0}$ belongs to the clean intersection calculus, see H\"ormander \cite{Ho4}. However,  as demonstrated in \cite[Lemma 11.1]{LOSU}, $\rho$ fails to be injective on the set $\mcl\cap C$ where
\beq
\mcl = \{(y, v, \eta, w; x, \xi) \in T^*\mcc_0 \times T^*M: \text{$\xi$ is light like}\}.
\eeq

Now let's consider the wave front set of the normal operator. The canonical relation for $L_{\mcc_0}^t$ is $C^t$, so by the calculus of wave front set (see e.g.\ \cite{Ho4}), we have
\beq
\begin{split}
\WF(L_{\mcc_0}^t \circ L_{\mcc_0}) &\subset (C\backslash \mcl)^t \circ (C\backslash \mcl)\cup (C\cap\mcl)^t \circ (C\cap\mcl)\\
&=\lap' \cup (C\cap \mcl)^t\circ (C\cap \mcl) \subset \lap' \cup \La'.
\end{split}
\eeq
Here we observed that 
\beq
\begin{gathered}
(C\cap \mcl)^t\circ (C\cap \mcl) = \{(x, \xi, z,  \zeta) \in (L^*M\backslash 0) \times (L^*M\backslash 0):\\
 \text{$(x, \xi)$ and $(z, \zeta)$ lie on a null bicharacteristics intersecting $L^*_\mcu M$}\}  \subset \La'.
 \end{gathered}
\eeq
To show that  $L_{\mcc_0}^t \circ L_{\mcc_0}$ actually belongs to the paired Lagrangian space $I^{p, l}(\lap, \La)$ and determine $p, l$, it suffices to show that the symbol belongs to the class of symbols of product type. For convenience, we shall work with $\chi L_{\mcc_0}$ for $\chi \in C_0^\infty(\mcu)$ and find the symbol of $(\chi L_{\mcc_0})^t \circ (\chi L_{\mcc_0})$. Here we can regard $\chi(y)$ as a function $\chi(y, v)$ defined on $\mcc_0$. Moreover, the analysis below works for any $\chi\in C_0^\infty(\mcc_0)$. 

For $f, h \in C_0^\infty(M; \sym^2)$, we compute
\beq
\begin{split}
(\chi L_{\mcc_0} f, \chi L_{\mcc_0} h)_{L^2(\mbr^3 \times S^2)} &= \int_{S^2}\int_\mbr \int_{\mbr^3} \int_\mbr \chi^2(y) f_{lm}(r, y+ rv) h_{jk}(s, y+ sv)\theta^j\theta^k\theta^l\theta^m ds dy dr dv\\
& = \int_{S^2}\int_\mbr \int_{\mbr^4} \chi^2(x\rq{}-x^0v) f_{lm}(r, x\rq{}+ (r-x^0)v) h_{jk}(x)\theta^j\theta^k\theta^l\theta^m dx dr dv,
\end{split}
\eeq
where we made the change of variable $x^0 = s, x\rq{} = y + sv$. We obtain that 
\beq
((\chi L_{\mcc_0})^t \circ (\chi L_{\mcc_0}) f)_{jk}(x) = \int_{S^2}\int_\mbr \chi^2(x\rq{}-x^0v)  f_{lm}(r, x\rq{}+ (r-x^0)v)  \theta^j\theta^k\theta^l\theta^m dr dv
\eeq
We can write this as an oscillatory integral using
\beq
f_{lm}(r, x\rq{}+ (r-x^0)v) = \frac{1}{(2\pi)^4} \int_{\mbr^4} \int_{\mbr^4} e^{i((r, x\rq{}+ (r-x^0)v) - z)\cdot\eta} f_{lm}(z) dzd\eta
\eeq
Therefore, we have
\beq
\begin{split}
&((\chi L_{\mcc_0})^t \circ (\chi L_{\mcc_0}) f)_{jk}(x)\\
= &\frac{1}{(2\pi)^4}\int_{S^2}\int_\mbr \int_{\mbr^4} \int_{\mbr^4} \chi^2(x\rq{}-x^0v) e^{ir(\eta^0 + v\cdot \eta\rq{})} e^{i\phi(x, z, v)}  \theta^j\theta^k\theta^l\theta^m f_{lm}(z)  dr dv dz d\eta,\\
 =&  \frac{1}{(2\pi)^4}\int_{S^2} \int_{\mbr^4} \int_{\mbr^4} \delta(\eta^0 + v\cdot \eta\rq{}) \chi^2(x\rq{}-x^0v) e^{i\phi(x, z, v)}  \theta^j\theta^k\theta^l\theta^m f_{lm}(z) dv dz d\eta
\end{split}
\eeq
where the phase function
\beq
\phi(x, z, v) = (-z^0, x\rq{}-x^0v - z\rq{})\cdot\eta = (x^0-z^0, x\rq{}-z\rq{})\cdot \eta,
\eeq
since the integrand is supported on $\eta^0 = -v\cdot \eta\rq{}$. Therefore, we can write 
\beq
\begin{split}
((\chi L_{\mcc_0})^t \circ (\chi L_{\mcc_0}) f)_{jk}(x) & =  \frac{1}{(2\pi)^4} \int_{\mbr^4} \int_{\mbr^4}   e^{i (x-z)\cdot \eta} a_{jklm}(x,  \eta)  f_{lm}(z) dz d\eta,
\end{split}
\eeq
where the symbol is given by 
\beq
a_{jklm}(x,  \eta) =  \int_{S^2}\delta(\eta^0 + v\cdot \eta\rq{}) \chi^2(x\rq{}-x^0v)  \theta^j\theta^k\theta^l\theta^m dv.
\eeq
The computation in \cite[Lemma 8.1]{LOSU} see also \cite[Prop.\ 11.4]{LOSU} showed that $a_{jklm}$ is a locally integrable function and the integral was explicitly evaluated, which we recall now. Consider the set $S^1_\eta = \{v\in S^2: \eta^0 + \eta\rq{} v = 0\}$.  If $\eta$ is time like, $S^1_\eta = \emptyset$. If $\eta$ is space-like, $S^1_\eta$ is a circle of radius $|\eta\rq{}|^{-1}(|\eta\rq{}|^2 - (\eta^0)^2)^{\ha}$. We have 
\beqq\label{eqsyma}
a_{jklm}(x, \eta) =\left\{ \begin{array}[c]{c}
  \dfrac{1}{(|\eta\rq{}|^2 - (\eta^0)^2)^{\ha}}\int_{S^1_\eta}   \chi^2(x\rq{}-x^0v) \theta^j\theta^k\theta^l\theta^m   dv, \text{ if $\eta$ is space-like},  \\ 
    \ \ \\
    0, \text{ otherwise.}\\
  \end{array}\right.
\eeqq
Now we see that on $\lap\backslash \La$, $a_{jklm}$ is a symbol of order $-1$ so that the normal operator is a pseudo-differential operator of order $-1$ microlocally restricted to $\lap\backslash \La$. This was obtained in \cite{LOSU}. Also, $a_{jklm}$  is singular at $\Sigma$ consisting of light-like vectors $\eta$. According to the discussion at the end of Section \ref{main}, the symbol $a$ belongs to the class $S^{m, m'}(M\times M; \lap, \La)$ with  $m' = 0$. Moreover, we have $m' = l-\ha = 0, p+l = -1$ and we solve that  $p =-\frac 32, l =  \ha.$ Therefore, $(\chi L_{\mcc_0})^t \circ (\chi L_{\mcc_0}) \in I^{-\frac 32, \ha}(\lap, \La)$.  Now we can apply \cite[Theorem 3.3]{GrU90} and a duality argument to obtain the Sobolev estimates of $\chi L_{\mcc_0}$. Thus we've proved 
\begin{theorem}\label{main0}
For any $\chi\in C_0^\infty(\mcu)$ (or $C_0^\infty(\mcc_0)$), the normal operator $(\chi L_{\mcc_0})^t \circ (\chi L_{\mcc_0}) \in I^{-\frac 32, \ha}(M \times M; \lap, \La)$. Also for $s\geq -\ha$,  $\chi L_{\mcc_0}: H^s(M)\rightarrow H_{loc}^{s+ \ha}(\mcc_0)$ is bounded. 
\end{theorem}
 
We remark that the Sobolev estimates can be seen from a more general result of Greenleaf and Seeger\footnote{The author thanks Prof.\ Greenleaf for pointing this out. }. In \cite{GrS}, the authors demonstrated (in Section 4) that in absence of conjugate points, the light ray transform on a general Lorentzian manifold is an FIO associated to a canonical relation where one projection is a submersion with folds, and the mapping properties of such operators are analyzed. We can apply \cite[Corollary 4.2]{GrS} to  $L_{\mcc_0}$ to obtain $L_{\mcc_0}: H_{comp}^s(M) \rightarrow H^{s+1/2+\eps}_{loc}(\mcc_0)$ for any $\eps>0$. However, one can check the proof of \cite[Theorem 1.1]{GrS} to conclude that the $\eps$ loss does not happen for $L_{\mcc_0}$ because the Hessian of the submersion with folds in this case is sign-definite.

\section{The parametrix construction}\label{secpara}
We prove Theorem \ref{main1}. Notice that since we shall consider the operator $L_{\mcc_0}$ acting on distributions in $\mcp(\mcu)$ so that $L_{\mcc_0}$ is injective, we actually have $L_{\mcc_0} = L_{\mcc}$ so we just need to consider the light ray transform $L_{\mcc}$. The analysis in Section \ref{secnorm} applies to this case by taking $\mcu = \mbr^3$ and $\chi = 1$. 
Notice that $\lap\backslash\Sigma$ has disjoint components 
\beq
\begin{gathered}
\lap^- = \{(x, \xi, x, -\xi) \in T^*M\backslash 0 \times T^*M\backslash 0: \text{$\xi$ is time like}\} = \Omega^-\times \Omega^-\\
\lap^+ = \{(x, \xi, x, -\xi) \in T^*M\backslash 0 \times T^*M\backslash 0: \text{$\xi$ is space like}\} = \Omega^+\times \Omega^+
\end{gathered}
\eeq 
so that $\lap = \lap^+\cup \lap^-\cup \Sigma$. We consider the set where the symbol $a$ in \eqref{eqsyma} (when $\chi = 1$) is elliptic. 
For $(x, \eta)\in \Omega^+$, consider $a(x, \eta): f_{lm}\rightarrow  a_{jklm}(x, \eta)f_{lm}$ as a linear map on $\sym^2_x$. Since $\chi = 1$ does not vanish identically on $S^1_\eta$, we know from \cite[Lemma 9.1]{LOSU} that  the kernel of the map is given by 
\beq
\mcn_x = \{c g(x) + \eta\otimes w + w\otimes \eta; c\in \mbr, w\in \mbr^4\}, \text{ for any } x \in M.
\eeq
Therefore,  $a |_{\Omega^+}$ is injective on $C^\infty(M; \sym^2)\backslash \mcn$. In particular,  one can find $b_{ijkl}(x, \eta)$ such that $b_{\alpha\beta jk}(x, \eta) a_{jklm}(x, \eta)|_{ \Omega^+} = \delta_{\alpha l}\delta_{\beta m}$ on $C^\infty (M; \sym^2)\backslash \mcn$. Since $a_{jklm}(x, \eta)$ is a symbol  of order $-1$ on $\Omega^+$, we can find $b_{jklm}(x, \eta)$ a symbol of order $1$ on $\Omega^+$.

Now we use the calculus of paired Lagrangian distribution to construct a parametrix for $L_{\mcc_0}$. The argument is quite standard as for elliptic pseudo-differential operators. We will use the symbol calculus \cite[Prop.\ 3.4]{GrU1} and the composition of $I^{p, l}$ for the flow out model \cite[Prop.\ 3.5]{GrU1}. These results can be found in \cite{AU,  DUV, GU} as well. 
First, we let $A_0 \in I^{\ha, \ha}(\lap, \La)$ be an operator with a symbol  $\sigma(A_0)(x, \eta) = b(x, \eta)$ on $\Omega^+$ and otherwise $0$. Then we have that acting on $\mcp(\mcu)$, 
\beq
A_0 \circ L_{\mcc}^t\circ L_{\mcc_0} - H 
\in I^{-\ha, -\ha}(\lap, \La) + I^{-\frac 32, \ha}(\lap, \La)
\eeq
where $H\in I^{- \frac{1}{2},  \ha}(\lap, \La)$. Also, we have $H-  \mathbbm{1}_{\Omega^+} \in \bigcap_{l} I^{-\ha, l}(\lap, \La) = I^{-\ha}(\La)$. Next, using the ellipticity of the symbol $a$ (on $\Omega^+$), we can follow the argument in \cite[Page 226-227]{GrU1} to get  $A\in I^{\ha, \ha}(\lap, \La)$ such that 
\beq
A  \circ  L_{\mcc}^t\circ  L_{\mcc_0} = \mathbbm{1}_{ \Omega^+} + B, \ \ B\in I^{-\ha}(\La)
\eeq
modulo a smoothing operator and acting on distributions in $\mcp(\mcu).$ This completes the proof of Theorem \ref{main1}.


\end{document}